\documentclass{jssc}

\def\textsubscript#1%
{$_{\text{#1}}$}

\def\cdd{\mbox{\boldmath$\cdot$}~}

\input{amssym.def}
\include{graphix}

\usepackage{color}
\usepackage{amsfonts,amssymb}
\usepackage{bm}
\usepackage{amsmath}
\usepackage{multirow}
\usepackage{algorithmic}
\usepackage{graphicx}
\usepackage{textcomp}
\usepackage{graphicx} 
\usepackage{epstopdf}

\newtheorem{mydef}{Definition}[section]
\newtheorem{Assumption}{Assumption}[section]

\usepackage{array}
\usepackage{enumerate}
\usepackage{booktabs}
\usepackage{algorithm}
\usepackage{algorithmic}
\usepackage{setspace}
\usepackage{float}
\usepackage{cases}
\usepackage[title]{appendix}

\makeatletter
\def\@oddfoot{\hfill}
\newcount\shumeicount
\def\setshumei#1#2#3{%
  \shumeicount=\count0
  \def\@oddhead{%
    \raise-5pt\hbox to0pt{\vrule width\hsize height 0pt depth 0.4pt\hss}\relax
    \ifnum \shumeicount=\count0
      \raise-7pt\hbox to0pt{\vrule width\hsize height 0pt depth 0.4pt\hss}\relax
      #1
    \else
      \ifodd\count0
        #2
      \else
        #3
       \fi
     \fi
  }%
}
\makeatother
\makeatletter
\def\@oddfoot{\hfill}
\newcount\shujiaocount
\def\setshujiao{%
  \shujiaocount=\count0
  \def\@oddfoot{%
      \ifodd\count0
      \else
      \fi
  }%
}
\makeatother
\def\title#1#2#3#4{{
  \vspace*{0.3cm}
  \begin{flushleft} \Large\bf #1\end{flushleft}
  \vspace*{-0.2cm}
      \begin{flushleft}
      \bf #2
      \end{flushleft}
      \footnotetext{\hspace{-6mm} #3\\ #4}}}

\def\dshm#1#2#3#4
{\setshumei{J Syst Sci Complex (#1) #2:
{\thepage--\pageref{LastPage}}\hfill}
            {\hfill {\small #3}\hfill\hbox to0pt{\hss\thepage}}
            {\hbox to0pt{\thepage\hss}\hfill {\small #4}\hfill
            }
            \setshujiao}
\def\drd#1#2
{{\vskip 1cm\small \begin{flushleft}
 #1 \\
 #2\\
\copyright The Editorial Office of JSSC \&  Springer-Verlag GmbH
Germany 2018
\end{flushleft}}}

\def\bar{\overline}
\def\epsilon{\varepsilon}

\begin{document}

\title{Algorithm design and approximation analysis  on distributed robust game$^*$}
{ \uppercase{XU } Gehui \cdd \uppercase{CHEN}
Guanpu \cdd \uppercase{QI} Hongsheng}
{\uppercase{XU } Gehui \cdd \uppercase{CHEN}
	Guanpu \cdd \uppercase{QI} Hongsheng (Corresponding author) \\
Key Laboratory of Systems and Control, Academy of Mathematics and Systems Science, Chinese Academy of Sciences, Beijing 100190,
and School of Mathematical Sciences, University of Chinese Academy of Sciences, Beijing 100049, China.   \\
Email: xghapple@amss.ac.cn; chengp@amss.ac.cn; qihongsh@amss.ac.cn   
   } 
{
{$^*$This work is supported partly by the National Key R\&D Program
of China under Grant 2018YFA0703800,  the Strategic Priority Research Program of Chinese Academy of Sciences under Grant No. XDA27000000, and the National Natural Science Foundation of China under Grants 61873262 and 61733018.}\\
{$^\diamond${\it This paper was recommended for publication by
Editor . }}}

\drd{DOI: }{Received: x x 20xx}{ / Revised: x x 20xx}


\dshm{20XX}{XX}{ALGORITHM DESIGN AND APPROXIMATION ANALYSIS  ON DISTRIBUTED ROBUST GAME}{\uppercase{XU } Gehui \cdd \uppercase{CHEN}
	Guanpu \cdd \uppercase{QI} Hongsheng}

\Abstract{We design a distributed algorithm to seek generalized Nash equilibria of a robust game with uncertain coupled constraints. Due to the  uncertainty of parameters in set constraints, we aim to find a generalized Nash equilibrium in the worst case. However,
	  it is challenging to obtain the exact equilibria directly because the parameters  are from general  convex sets, which may not have analytic expressions or are endowed with  high-dimensional nonlinearities.
	To solve this problem, we first  approximate parameter sets with inscribed polyhedrons, and  transform the approximate problem in the worst case  into an extended certain game with resource allocation constraints by robust optimization. Then we propose a distributed algorithm for this certain game and prove that an equilibrium obtained from the algorithm induces an $ \epsilon $-generalized Nash equilibrium of the original game, followed by convergence analysis.
	Moreover, resorting to the metric spaces and the
	analysis on nonlinear perturbed systems,
	we estimate the approximation accuracy related to $ \epsilon $ and point out the factors
	influencing the accuracy of $ \epsilon $.}        




\section{Introduction}


Multi-agent systems involving a non-cooperative setting
have attracted extensive research and applications in many
fields, such as telecommunication power allocation and cloud
computation  \cite{ardagna2012generalized,pang2008distributed}.   Due to some shared resources between players, such as communication bandwidth and network energy, coupled constraints are frequently considered in non-cooperative games.
As a reasonable solution, a generalized Nash equilibrium (GNE) can be regarded as  defined as a set of strategies that satisfies the local and coupled constraints, in which no
player can profit from unilaterally deviating from its own strategy.  Significant
theoretic and algorithmic achievement of GNE seeking have
been done, referring to \cite{facchinei2007finite, fischer2014generalized}.


Recently, seeking equilibria in a distributed manner has
become an emerging research topic, where players obtain the
Nash equilibrium (NE) or GNE by making decisions with
local information and communicating through networks. Various distributed
algorithms have been proposed for  GNE seeking, such as   asymmetric projection algorithms \cite{paccagnan2016distributed},  projected dynamics based on non-smooth tracking dynamics \cite{liang2017distributed}, and forward–backward operator splitting method \cite{yi2019operator} with extended to fully distributed games  \cite{belgioioso2020distributed}.

However, considering the impact of the inevitable uncertainties
in practical games, it is often difficult to obtain the
exact GNE  directly in practice. One
way to handle uncertainties is to utilize robust optimization \cite{bertsimas2011theory}, which addresses the robust counterpart of an optimization model with uncertain data/parameters. By employing the robust optimization
approaches to deal with the uncertainties in games,
the concept of robust game was first proposed in \cite{aghassi2006robust}. Hereupon, the works themed on robust game have been applied in various scenarios, such as  human decision-making models in security setting,  defensive resource allocation  in homeland security, downlink power control problem with  interfering channel information, and electric vehicle charging problem under demand uncertainty \cite{pita2012robust,nikoofal2012robust,zhu2014downlink,yang2015noncooperative}.

Nevertheless, the analysis of robust games with coupled constraints is less. Most of the previous works focused
on the uncertainties in payoff functions or strategy variables,
and very few studied the uncertainties in the parameters of the accompanied constraints.  In addition, considering that coupled constraints often occur in  actual games, distributed GNE seeking in robust games deserve further investigation.
More recently, \cite{chen2021distributed} studied a robust game with parameters
uncertainty in coupled constraints, where an approximation
method was proposed to find an $ \epsilon $-GNE of the original game
in the worst case, but the estimation of $ \epsilon $ was not considered.
As the approximation focuses on the parameter sets while $ \epsilon $ is
affected by the feasible sets, it is hard to
construct the relationship between the approximation accuracy
and  $ \epsilon $. Furthermore, the difficulty of solving
the problem increases due to estimating $ \epsilon $ in a distributed setting. Therefore,
the  distributed robust game  with general uncertainty is hard
to be analyzed using the existing methods.


In this work, we study distributed  GNE seeking of a robust
game with general uncertainties, where the parameters in
coupled constraints are from general uncertain convex sets, which is more generalized than the previous works without uncertainty
in constraints \cite{paccagnan2016distributed,liang2017distributed,gadjov2018passivity}, or  restricted to special
structure \cite{zeng2018distributed11,wang2014generalized}.
Due to the complexity of uncertainty modeling, the parameter
sets may not be equipped with exact analytic expressions or are endowed with high-dimensional nonlinearities,
which makes it hard to obtain the exact equilibria directly. To
solve this problem, we approximate uncertain  parameter sets  with inscribed polyhedrons and transform  the approximate problem in
the worst case into an extended certain game model with resource allocation constraints by robust optimization. Then we propose a distributed continuous-time
algorithm for seeking a GNE of the  certain game,
followed by the convergence analysis.  The proposed algorithm has lower dimensions than \cite{chen2021distributed},  and avoids discontinuities caused by  tangent cones in \cite{chen2021distributed,bianchi2021continuous}.
Moreover, by virtue of  metric spaces and perturbed systems, an equilibrium obtained from the algorithm is proved to be an $ \epsilon $-GNE of the original game, and an upper bound of the approximation accuracy related to $ \epsilon $ is given.
The remainder is organized as follows. Section \ref{s2} provides
notations and preliminary knowledge, while Section \ref{s3} formulates
a  distributed robust game with parameter uncertainties
in coupled constraints. Then Section \ref{s4} provides a distributed
algorithm based on a resource allocation
problem after a proper approximation and gives the convergence
analysis. Section \ref{s6} shows that the equilibria of the
designed algorithm are $ \epsilon $-GNE of the original problem in the worst case  and
obtains an upper bound of the value $ \epsilon $, and Section \ref{s7} presents
numerical examples for illustration of the proposed algorithm.
Finally, Section \ref{s8} concludes the paper.

\section{Preliminaries}\label{s2}
In this section, we introduce some basic notations  and  preliminary knowledge.

Denote $ \mathbb{R}^{n} $ (or $ \mathbb{R}^{m\times n} $) as the set of $ n $-dimensional (or $ m $-by-$ n $) real column vectors (or real matrices), and $ I_{n} $ as the $ n\times n $ identity matrix. Let $ \boldsymbol{1}_{n} $(or $ \boldsymbol{0}_{n} $) be the $ n $-dimensional column vector with all elements of $ {1} $ (or $ {0} $).   For a column vector $ x \in \mathbb{R}^{n}$, $ x^{\mathrm{T}} $ denotes
its transpose.  Take $ \operatorname{col}\{x_{1},\cdots,x_{n}\} $$  = (x^{\mathrm{T}}_{1}, \cdots,x^{\mathrm{T}}_{n})^{\mathrm{T}} $ as the stacked column vector obtained from column vectors $ x_{1},\cdots,x_{N} $, $\|\cdot\|$ as the Euclidean norm, and $ \operatorname{relint}(D) $ as the relative interior of the set $ D $. Denote $ \operatorname{ker}(M) $ as the kernel of the matrix $ M $, $ \operatorname{Im}(M) $ as the image space of the matrix $ M $ and $ \operatorname{span}(x) $ as the spanning subspace by vector $ x $.
Denote $ \mathbf{E}_{v}(c)\subseteq \mathbb{R}^{n} $ as an ellipsoid that
$$
\sum_{i=1}^{n}\frac{(x_{i}-c_{i})^{2}}{v_{i}^{2}}\leq1,\;
$$
with the center at point $ c\triangleq(c_{1},\cdots,c_{n}) $ and the semiaxis  $v\triangleq(v_{1},\cdots,v_{n})$.

A set $ \Omega \subseteq \mathbb{R}^{n} $ is convex if $ \omega x_{1}+(1-\omega)x_{2} \in \Omega$ for  any
$ x_{1}, x_{2}\in \Omega $ and $ 0\leq\omega\leq 1   $.
For a closed convex set $ \Omega $, the projection map $ \Pi_{\Omega}: \mathbb{R}^{n} \rightarrow \Omega $ is defined as
$$
\Pi_{\Omega}(x) \triangleq \underset{y \in \Omega}{\operatorname{argmin}}\|x-y\|.
$$
Especially, denote $ [x]^{+}\triangleq \Pi_{\mathbb{R}^{n}_{+}}(x) $ for convenience.

A mapping $ F : \mathbb{R}^{n} \rightarrow \mathbb{R}^{n} $ is said to be   monotone  (strictly monotone) on a set $ K $ if
$$
(F(x)-F(y))^{\mathrm{T}}(x-y) \geq 0\, (>0), \quad \forall x, y \in K, x\neq  y.
$$

Given a set $ K \subseteq \mathbb{R}^{n} $ and a map $ F: K \rightarrow \mathbb{R}^{n} $, the variational inequality problem $\mathrm{VI}(K, F)$ is defined to find a vector $ x^{*}\in K  $ such that $$
\left(y-x^{*}\right)^{\mathrm{T}} F\left(x^{*}\right) \geq 0, \quad \forall y \in K,
$$
whose  solution  is denoted by $\mathrm{SOL}(K, F)$.  When $ K $ is closed and convex, the  solution of $\mathrm{VI}(K, F)$  can be equivalently reformulated via projection as
$$
x \in \operatorname{SOL}(K, F) \Leftrightarrow x=\Pi_{K}(x-  F(x)).
$$
Moreover, if $ K $ is compact, then $\mathrm{SOL}(K, F)$ is nonempty and compact. If $ K $ is closed and $ F(x) $ is strictly monotone, then $\mathrm{VI}(K, F)$ has at most one solution \cite[Proposition 1.5.8, Corollary 2.2.5, and Theorem 2.3.3]{facchinei2007finite}.

Take $ X,\, Z\subseteq  \mathbb{R}^{n} $ as two non-empty sets. For $ y \in \mathbb{R}^{n}  $,  denote $ \operatorname{dist}(y,Z) $ as the distance between $ y $ and $ Z $, i.e.,
$$ \operatorname{dist}(y,Z) = \inf \limits_{z\in Z}\|y- z\|.  $$
Define the Hausdorff metric of $ X,Z\subseteq  \mathbb{R}^{n} $ by $$ H(X,Z) = \max\{\sup\limits_{x\in X}\operatorname{dist}(x,Z),\sup\limits_{z\in Z}\operatorname{dist}(z,X) \}.$$The Hausdorff metric integrates all compact sets into a metric space.

Let $ \mathcal{X} $ and $\mathcal{Y} $ be $ m $-dimensional subspaces of $ \mathbb{R}^{n} $, respectively. The canonical angles
between them are defined to be
$$
\vartheta_{i}(\mathcal{X}, \mathcal{Y})=\arccos \sigma_{m-i+1}\left(X^{\mathrm{T}} Y\right), \; i=1,2, \ldots, m,
$$
where $ X $ and $ Y $ are matrices whose columns form orthonormal bases of $ \mathcal{X} $ and $\mathcal{Y} $, and $ \sigma_i(X^{\mathrm{T}} Y) $,
$ i=1,2, \ldots, m $, are decreasingly ordered singular values of $ X^{\mathrm{T}} Y $. Denote the
canonical angles between $ \mathcal{X} $ and $ \mathcal{Y} $ by
$
\vartheta(\mathcal{X}, \mathcal{Y})\triangleq (\vartheta_{1}(\mathcal{X}, \mathcal{Y}) ,...,\vartheta_{m}(\mathcal{X}, \mathcal{Y})).
$
The following lemma reveals  the metric about canonical angles between $ \mathcal{X} $ and $ \mathcal{Y} $ \cite{qiu2005unitarily}, \cite{zhang2007angular}.
\begin{lemma}\label{l5}
	Let $ \varrho : \mathbb{R}^{m}\rightarrow \mathbb{R} $ be a symmetric gauge function.  Define $ \psi : \mathbb{R}^{m}\times\mathbb{R}^{m} \rightarrow \mathbb{R} $ of $ \mathcal{X} $ and $\mathcal{Y} $ by
	$$ \psi(\mathcal{X},\mathcal{Y})=\varrho(\vartheta(\mathcal{X}, \mathcal{Y})). $$
	Then $ \psi $  is called  an angular metric. Moreover,
	let $ \mathcal{X}_{\perp} $ and $ \mathcal{Y}_{\perp} $ be the orthogonal complements of $ \mathcal{X} $ and $ \mathcal{Y} $, respectively. The nonzero canonical angles between $ \mathcal{X} $ and $ \mathcal{Y} $ are the same as those of $ \mathcal{X}_{\perp} $ and $ \mathcal{Y}_{\perp} $, which means that $ \psi(\mathcal{X},\mathcal{Y}) =\psi(\mathcal{X}_{\perp},\mathcal{Y}_{\perp})$.
\end{lemma}

Consider a class of comparison functions.
A continuous function  $ \alpha : [0, a) \rightarrow [0,\infty)  $  is said to belong to class $ \mathcal{K} $ if it is strictly increasing and $ \alpha(0)=0 $. It is said to belong to class  $ \mathcal{K}_{\infty} $ if $ a=\infty $ and $ \alpha(r)\rightarrow  \infty$ as $ r\rightarrow \infty $.

Moreover,  the information sharing among the players can be described by a graph $ \mathcal{G} = (\mathcal{I},\mathcal{E}) $, with the node set  $ \mathcal{I} =\{1,2,\cdots,N\} $ and the edge set $ \mathcal{E} $. $ A = [a_{i j}] \in \mathbb{R}^{n\times n} $ is the adjacency matrix of $ \mathcal{G} $ such that if $ (j, i) \in \mathcal{E}  $, then $ a_{i j}> 0 $,   which
means that $ i $ can obtain the information from $ j $ and
$j$ belongs to  $ i $'s neighbor set;  $ a_{i j}= 0 $ otherwise. $ \mathcal{G} $ is said to be undirected if $ (j, i) \in \mathcal{E} \Leftrightarrow (i, j) \in \mathcal{E} $, and  $ \mathcal{G} $ is  to be connected if any two nodes in $ \mathcal{I} $ are connected by a path.
The Laplacian matrix is $ L = \Delta-A $, where $\Delta=\operatorname{diag}\left\{d_{1}, \ldots, d_{N}\right\} \in \mathbb{R}^{N \times N}$ with $ d_{i}=\sum_{j=1}^{N} a_{i j} $.  When $ \mathcal{G} $ is an undirected connected  graph, $ 0 $ is a
simple eigenvalue of Laplacian $ L $ with the eigenspace $ \{a \boldsymbol{1}_{n} |a \in
\mathbb{R}\} $, and $ L \boldsymbol{1}_{n} = \boldsymbol{0}_{n} $, while all other eigenvalues are positive.

\section{Problem Formulation}\label{s3}
Consider an $ N $-player  game with a global coupled constraint  as follows.  For  $i \in \mathcal{I}\triangleq \{1,\cdots,N\}$,  player $i$ has an action variable $x_{i}$ in a local action  set  $\Theta_{i}\subseteq \mathbb{R}^{n} $. Denote $  \boldsymbol{\Theta} =\prod_{i=1}^{N} \Theta_{i}\subseteq \mathbb{R}^{nN} $, $ \boldsymbol{x}\triangleq \operatorname{col}\{x_{1}, . . . ,x_{N}\} \in \boldsymbol{\Theta} $ as the action profile for all players, and $ \boldsymbol{x}_{-i}\triangleq \operatorname{col}\{x_{1}, . . . ,x_{i-1}, x_{i+1}, . . . ,x_{N}\}$ as the action profile for all players except player $ i $.  The  cost function for player $ i $ is $  J_{i}(x_{i},\boldsymbol{x}_{-i}) : \mathbb{R}^{nN}\rightarrow \mathbb{R} $.

Moreover,
there exists a coupled inequality constraint shared by all players. Denote $ \boldsymbol{K}\subseteq \mathbb{R}^{Nn}$ as the   set for this coupled constraint. Considering  that the parameters in
constraints are given in general uncertain convex sets, the action profile $ \boldsymbol{x} $ needs to satisfy
$$
\boldsymbol{x}\!\in\boldsymbol{K} \!\triangleq\!\left\{\boldsymbol{x} \!\in \mathbb{R}^{Nn} \,\Big|\,\sum_{i=1}^{N} \omega_{i}^{\mathrm{T}} x_{i} \!\leq\! b, \; \omega_{i}\! \in\! \mathcal{M}_{i} \!\subseteq\! \mathbb{R}^{n},\; \forall i \in \mathcal{I}\right\},
$$
where  $ \mathcal{M}_{i} $ is convex and compact. 
For any $ \omega_{i} \in \mathcal{M}_{i} $, the inequality constraint must be satisfied.
Denote the feasible action set of this game
by $ \boldsymbol{\mathcal{X}} \triangleq \boldsymbol{K}\bigcap \boldsymbol{\Theta}. $ Then, the feasible set of player $ i $ is
$$
\mathcal{X}_{i}(\boldsymbol{x}_{-i})\triangleq \left\{x_{i} \in \Theta_{i} \,\Big|\, \omega_{i}^{\mathrm{T}}  x_{i} \leq b-\!\!\sum_{j \neq i, j \in \mathcal{I}} \omega_{j}^{\mathrm{T}}  x_{j}, \omega_{i} \in \mathcal{M}_{i} \right\}.
$$

To sum up, given $ \boldsymbol{x}_{-i} $, the $ i $th player aims to solve

\begin{equation}\label{f5}
	\begin{aligned}
		& \min _{x_{i} \in \mathbb{R}^{n}} J_{i}\left(x_{i}, \boldsymbol{x}_{-i}\right)
		\text { s.t. } x_{i} \in \mathcal{X}_{i}(\boldsymbol{x}_{-i}).
	\end{aligned}
\end{equation}

\begin{mydef}[$ \epsilon $-generalized Nash equilibrium]\label{d3}
	A profile $ \boldsymbol{x}^{*} $ is said to be an $ \epsilon $-generalized Nash equilibrium of game (\ref{f5})  if
	\begin{equation}\label{e6}
		J_{i}\left(x_{i}^{*}, \boldsymbol{x}_{-i}^{*}\right) \leq J_{i}\left(x_{i}, \boldsymbol{x}_{-i}^{*}\right)+\epsilon,\;\forall i \in \mathcal{I},\;\forall x_{i}\in \mathcal{X}_{i}(\boldsymbol{x}_{-i}),\;
	\end{equation}
	with a positive constant $ \epsilon $. Particularly, $ \boldsymbol{x}^{*} $  is said to be a GNE when $ \epsilon= 0 $.
\end{mydef}

The main task of this paper is to design a distributed dynamics
for seeking a GNE of the robust game (\ref{f5}), where each
player can only access  its local payoff function and feasible
decision set under a multi-agent network. The $ i $th player may only know  $ \omega_i^{\mathrm{T}}x_{i} $ and the parameter uncertainty set $ \mathcal{M}_{i} $,  rather than $ \sum_{i=1}^{N} \omega_i^{\mathrm{T}}x_{i} $.  To fulfill the cooperations between players
for solving (\ref{f5}), the players have to share their local information
through a network $ \mathcal{G} $.  On the other hand, restricted by the
uncertainty of $ \omega_{i} $, we aim to find a GNE of (\ref{f5}) in the worst
case, i.e., a GNE satisfies all possible constraints, which is defined as $$ \boldsymbol{x}^{*}\in  \left\{\boldsymbol{x} \in \boldsymbol{\Theta}\, \Big| \,\sum_{i=1}^{N} \max_{\omega_{i} \in 	\mathcal{M}_{i}} \omega_{i}^{\mathrm{T}} x_{i} \leq b \right\}.$$
%
However, it is very difficult to solve the worst-case solution directly, because the
	challenge  comes from the fact that $ \omega_{i} $  is arbitrarily selected from a general uncertain convex set $ \mathcal{M}_{i} $, which may
	be endowed with high-dimensional nonlinearities or have no
	analytical expression.
Therefore, we consider finding an $ \epsilon $-GNE of game (\ref{f5}) in the
worst case with a practical approximation, and analyze the
approximation accuracy related to $ \epsilon $, which overcomes the
difficulty of estimating $ \epsilon $ in \cite{chen2021distributed}.
%

\begin{remark}
	In our distributed game, the decision variable $ x_j $ can be observable by the $ i $th player, if $J_{i}(x_{i}, \boldsymbol{x}_{-i})$ depends explicitly on $x_{j}$, for any $j\in \mathcal{I}$.
	Thus, player $ i $ can get its local gradient by observing the decisions influencing $J_{i}(x_{i}, \boldsymbol{x}_{-i})$.
	This observation	model has also been adopted in \cite{chen2021distributed,yi2019operator}. On the other hand, there have also been methods for distributed GNE seeking when each player cannot observe the full decisions that its cost function depends on, referring to \cite{bianchi2021continuous,zhu2021generalized}. Here, we do not consider this circumstance, where this simplification does not affect the focus of our research.
\end{remark}

The following assumptions  are associated with game (\ref{f5}).
\begin{Assumption}\label{a3}
	$\ $
	\begin{itemize}
		\item For $ i \in \mathcal{I} $, $ \Theta_{i} $ is compact and convex.	Besides, there exists $ \boldsymbol{x}\in
		\text{relint}(\boldsymbol{\Theta}) $ such that
		$
		\sum_{i=1}^{N} \omega_{j}^{\mathrm{T}} x_{j} < b, \quad \omega_{j} \in \mathcal{M}_{j} \subseteq \mathbb{R}^{n}, \quad \forall j \in \mathcal{I}.
		$
		\item For $i\in \mathcal{I}$,  $ J_{i}( \boldsymbol{x}) $ is  Lipschitz continuous in $ \boldsymbol{x} $, while $ J_{i}( \boldsymbol{x}) $ is  continuously differentiable in $ x_i $.
		Moreover, the pseudo-gradient $F(\boldsymbol{x}) \triangleq \operatorname{col}\left\{\nabla_{x_{1}} J_{1}\left(\cdot, \boldsymbol{x}_{-1}\right), \ldots, \nabla_{x_{N}} J_{N}\left(\cdot, \boldsymbol{x}_{-N}\right)\right\}$ is strictly monotone in $ \boldsymbol{x} $.
		\item The undirected graph  $ \mathcal{G} $ is   connected.
	\end{itemize}
\end{Assumption}

By Assumption \ref{a3}, it is clear that Slater’s condition is
satisfied \cite{belgioioso2020distributed, zeng2019generalized}. Besides, compared with \cite{chen2021distributed,liang2020distributed}, the  map $ F $ is assumed to be strictly monotone rather than strongly
monotone.
\section{Algorithm Design}\label{s4}
In this section, we approximate the parameter uncertainty  sets of game (\ref{f5}) in a proper way and propose a distributed algorithm to find the worst-case solution with the
uncertainty in the approximate game.

 One of the most common tools for approximating convex sets is by inscribed polyhedrons \cite{kamenev1992class,kamenev2003self}.	
	 Recalling the definition of inscribed polyhedrons, it is a  polyhedron with all its vertices on the boundary of the convex set. And it is essentially   enclosed by a series of hyperplanes. Denote $  \boldsymbol{\mathcal{M}} =\prod_{i=1}^{N} \mathcal{M}_{i} $ and  $  \boldsymbol{\mathcal{P}}_{v} =\prod_{i=1}^{N} \mathcal{P}_{v_{i}}^{i} $.  Take  $ \mathcal{P}_{v_{i}}^{i} $  as an  inscribed polyhedron of $ \mathcal{M}_{i}$ with $ v_{i} $ vertices, it can be   expressed as
\begin{equation}\label{bbs}
	\mathcal{P}_{v_{i}}^{i}=\left\{\omega_{i} \in \mathbb{R}^{n}\,\Big|\, A_{i} \omega_{i} \leq d_{i}\right\}.
\end{equation}
Here, for $ i \in \mathcal{I} $, $ A_{i} \in \mathbb{R}^{q_{i}\times n}$  are normal vectors of the hyperplanes with normalized rows. They determine the directions of these hyperplanes. $ q_{i} $ is the number of hyperplanes, and $ d_{i} $ are the distances from the origin point to the hyperplanes.

\begin{remark}
Here we choose polyhedrons for the approximation because they can be explicitly expressed by linear inequalities, which provide
simple mathematical derivation and make the distributed algorithms concise. Furthermore, although  the analytical expressions of convex sets with high-dimensional nonlinearities are hard to solve directly,  in some situations one can sample exactly a few points on the boundary of the convex set, which naturally form an inscribed polyhedron.  This is another important reason for choosing inscribed polyhedrons.
\end{remark}	

 With the help of the approximation by inscribed polyhedrons, the
coupled constraint of (\ref{f5}) in the worst case becomes


\begin{equation}\label{ddd}
	\sum_{i=1}^{N} \max_{\omega_{i} \in 	\mathcal{P}_{v_{i}}^{i}} \omega_{i}^{\mathrm{T}} x_{i} \leq b .
\end{equation}
Then we can explicitly investigate the worst-case solution with
uncertainty based on robust optimization \cite{bertsimas2011theory} and robust game \cite[Theorem 1]{chen2021distributed}. Specifically, by introducing a dual variable  $ \sigma_{i}\in\mathbb{R}^{q_{i}}_{+} $,  (\ref{ddd}) can be equivalently transformed into
%
%
\begin{equation}\label{f4}
	\begin{array}{l}
		\sum^{N}_{j=1}d_{j}^\mathrm{T}\sigma_{j}\leq b,\quad A_{j}^\mathrm{T}\sigma_{j}-x_{j}= \boldsymbol{0}_{n},\quad \forall j\in \mathcal{I}.
	\end{array}
\end{equation}



Moreover, denote $ z_{i}=\operatorname{col}\{x_{i}, \sigma_{i}\} \in \mathbb{R}^{n+q_{i}}$,  $ 	B_{i}=\left[\mathbf{0}_{n}^{\mathrm{T}}, d_{i}^{ \mathrm{T}}\right] \in \mathbb{R}^{1 \times\left(n+q_{i}\right)}, $ and $	C_{i}=\left[-I_{n}, A_{i}^{ \mathrm{T}}\right] \in \mathbb{R}^{n \times\left(n+q_{i}\right)}$. Define $\Phi_{i}=\Theta_{i}\times \mathbb{R}^{q_{i}} $,
\begin{equation}\label{fff}
	\Omega_{i}=\Phi_{i}\cap\{C_{i}z_{i}=\boldsymbol{0}_{n}\},
\end{equation}
$ \boldsymbol{z}_{-i} $ as
all the vectors except $ z_{i} $, $ \boldsymbol{z}\triangleq \operatorname{col}\{z_{1}, . . . ,z_{N}\} \in \mathbb{R}^{nN+q} $, where $ q=\sum^{N}_{i=1}q_i $.
With these notations, game (\ref{f5}) with approximation is therefore converted into an extended certain game model with
resource allocation constraints, that is,
\begin{equation}\label{f7}
	\begin{array}{l}
		\min \limits_{z_{i} \in \Omega_{i}} \widehat{J}_{i}\left(z_{i}, \boldsymbol{z}_{-i}\right) \\
		\text { s.t. }  \sum_{j=1}^{N} B_{j} z_{j} \leq b, \quad \forall j \in \mathcal{I},
	\end{array}
\end{equation}
where  $ \widehat{J}_{i}\left(z_{i}, \boldsymbol{z}_{-i}\right)=J_{i}\left(x_{i}, \boldsymbol{x}_{-i}\right)  $.


%


Denote the pseudo-gradient of (\ref{f7}) by
$$
\boldsymbol{g}(\boldsymbol{z}) \triangleq \operatorname{col}\left\{g_{1}\left(z_{1}, \boldsymbol{z}_{-1}\right), \ldots, g_{N}\left(z_{N}, \boldsymbol{z}_{-N}\right)\right\} \in \mathbb{R}^{n N+q},
$$
where  $g_{i}\left(z_{i}, \boldsymbol{z}_{-i}\right) \triangleq \operatorname{col}\left\{\nabla_{x_{i}} J_{i}\left(\cdot, \boldsymbol{x}_{-i}\right), \mathbf{0}_{q_{i}}\right\} \in \mathbb{R}^{n+q_{i}}$.
Take $\boldsymbol{B}=\operatorname{Diag}\left(B_{1}, \ldots, B_{N}\right) \in \mathbb{R}^{N \times(n N+q)}$, $\boldsymbol{b}=\operatorname{col}\left\{b_{1}, \ldots, b_{N}\right\} \in \mathbb{R}^{N} $ with $ \sum_{i=1}^{N} b_{i}=b $.
Then the  feasible set  of  player $ i $ in (\ref{f7}) is defined as $$
\Xi_{i}\left(z_{-i}\right) \triangleq\left\{z_{i} \in \Omega_{i} \,\Big|\,\sum_{j=1}^{N} B_{j} z_{j} \leq b\right\}.
$$

Let $  \boldsymbol{\Xi} =\prod_{i=1}^{N} \Xi_{i} $, $  \boldsymbol{\Omega} =\prod_{i=1}^{N} \Omega_{i} $ and $  \boldsymbol{\Phi} =\prod_{i=1}^{N} \Phi_{i} $. Referring to \cite[Proposition 1.4.2]{facchinei2007finite} and \cite{facchinei2010generalized}, a strategy profile $ \boldsymbol{z}^{*} $ is said to be a variational equilibrium, or variational GNE,  if $ \boldsymbol{z}^{*} \in \mathrm{SOL}(\boldsymbol{\Xi}, \boldsymbol{g}(\boldsymbol{z})) $. Moreover, for a variational GNE of  game (\ref{f7}),
$ \boldsymbol{z}^{*}$ together with  multiplier $ \boldsymbol{\lambda}^{*}$ satisfy the following first order conditions,
\begin{subequations}\label{aaf}
	\begin{align}
		&	\boldsymbol{0}_{nN}  \in\boldsymbol{g}(\boldsymbol{z}^{*})+\boldsymbol{B}^{\mathrm{T}} \boldsymbol{\lambda}^{*}+\mathcal{N}_{\boldsymbol{\Omega}}(\boldsymbol{z}^{*}), \\
		&	0 \leq-(\boldsymbol{B} \boldsymbol{z}^{*}-\boldsymbol{b})^{\mathrm{T}} \cdot \boldsymbol{1}_{N}, \quad 0=(\boldsymbol{B} \boldsymbol{z}^{*}-\boldsymbol{b})^{\mathrm{T}} \boldsymbol{\lambda}^{*},\\
		&	\boldsymbol{0}_{N}   =L\boldsymbol{\lambda}^{*},
	\end{align}
\end{subequations}
where multiplier $ \boldsymbol{\lambda}^{*}= \operatorname{col}\{\lambda_{1}^{*},\cdots, \lambda_{N}^{*}\}\in \mathbb{R}_{+}^{N}$, and
$ L $ is   the Laplacian matrix of network $ \mathcal{G} $.



By solving the first order conditions (\ref{aaf}) of the variational inequality $ \mathrm{VI}(\boldsymbol{\Xi}, \boldsymbol{g}(\boldsymbol{z})) $, we derive a variational GNE of game (\ref{f7}), which can be regarded as a GNE with
	equal multipliers, i.e, $ \lambda_{i}^{*}= \lambda_{j}^{*} $,  $ \forall i,j \in \mathcal{I} $.

Furthermore, by employing an additional variable $ \boldsymbol{\zeta}= \operatorname{col}\{\zeta_{1},\cdots, \zeta_{N}\} \in \mathbb{R}^{N}$, we propose a distributed algorithm for solutions to (\ref{aaf}) of  approximate game (\ref{f7}).

\begin{algorithm}[H]
	\small
	\caption{for each  $ i \in \mathcal{I}  $}
	\label{a1}
	\vspace{0.1cm}
	
	\textbf{Initialization}:
	\vspace{-0.2cm}
	\begin{flalign*}
		& z_{i}(0)\in \Omega_{i},\;\lambda_{i}(0)\in \mathbb{R}_{+},\; \zeta_{i}(0)\in\mathbb{R}.
		&
	\end{flalign*}
	\textbf{Dynamics renewal}:	
	\begin{flalign*}
		&\dot{z}_{i}=\Pi_{\Omega_{i}}\left(z_{i}-g_{i}\left(z_{i}, \boldsymbol{z}_{-i}\right)-B_{i}^{\mathrm{T}} \lambda_{i}\right)-z_{i}, \\
		&\dot{\lambda}_{i}=\left[\lambda_{i}\!+B_{i} z_{i}-b_{i}-\sum_{j=1}^{N} a_{i j}\left(\lambda_{i}-\lambda_{j}\right)-\sum_{j=1}^{N}\! a_{i j}\left(\zeta_{i}-\zeta_{j}\right)\right]^{+}-\lambda_{i}, \\
		&\dot{\zeta}_{i}=\sum_{j=1}^{N} a_{i j}\left(\lambda_{i}-\lambda_{j}\right),
	\end{flalign*}
	where $ a_{ij} $ is the $ (i,j) $th element of the adjacency matrix.
\end{algorithm}
%

Equivalently,
a compact form of Algorithm 1 can be written as
\begin{equation}\label{e10}
	\left\{\begin{array}{ll}
		\dot{\boldsymbol{z}}=\Pi_{\boldsymbol{\Omega}}\left(\boldsymbol{z}-\boldsymbol{g}(\boldsymbol{z})-\boldsymbol{B}^{\mathrm{T}}\boldsymbol{\lambda} \right)-\boldsymbol{z}, & \boldsymbol{z}(0) \in \boldsymbol{\Omega}, \\
		\dot{\boldsymbol{\lambda}}=\left[\boldsymbol{\lambda}+\boldsymbol{B}\boldsymbol{z}-\boldsymbol{b}-L\boldsymbol{\lambda}-L\boldsymbol{\zeta}\right]^{+}-\boldsymbol{\lambda}, & \boldsymbol{\lambda}(0)=\mathbb{R}_{+}^{N},\\
		\dot{\boldsymbol{\zeta}}=L\boldsymbol{\lambda},& \boldsymbol{\zeta}(0)\in\mathbb{R}^{N}.\\
	\end{array}\right.
\end{equation}

In Algorithm 1, the $ i $th player calculates the local decision
variable $ z_{i}\in \Omega_{i} $  based on projected gradient play
dynamics. The local variable $ \lambda_{i}\in \mathbb{R}_{+} $ is to  estimate a
dual variable associated with the coupled constraints, while
the local auxiliary variable $ \zeta_{i}\in \mathbb{R} $  is calculated for the consensus of  $ \lambda_{i} $.
\begin{remark}
	Compared with the algorithm in \cite{chen2021distributed}, dynamics (\ref{e10}) is with lower
	dimensions. Meanwhile, (\ref{e10}) adopts the projection
	operation to deal with local feasible constraints, which
	avoids the discontinuous dynamics caused by tangent cones
	\cite{chen2021distributed,bianchi2021continuous}.
\end{remark}
The following lemma shows the equivalence between an
equilibrium of algorithm (\ref{e10}) and a solution to $ \mathrm{VI}(\boldsymbol{\Xi}, \boldsymbol{g}(\boldsymbol{z})) $
satisfying (\ref{aaf}).
\begin{lemma}\label{l3}
	Under Assumption \ref{a3}, consider the game (\ref{f7}). If
	$ \operatorname{col}\{\boldsymbol{z}^{*}, \boldsymbol{\lambda}^{*}, \boldsymbol{\zeta}^{*}\}  $ is an equilibrium point of (\ref{e10}), then $ \boldsymbol{z}^{*} $ is a
	variational GNE of (\ref{f7}).  Conversely, if $ \boldsymbol{z}^{*} $ is a variational
	GNE of (\ref{f7}),  there exists $( \boldsymbol{\lambda}^{*}, \boldsymbol{\zeta}^{*})\in \mathbb{R}^{N}_{+} \times\mathbb{R}^{N} $
	such that $(\boldsymbol{z}^{*}, \boldsymbol{\lambda}^{*}, \boldsymbol{\zeta}^{*}) $ is an equilibrium point of (\ref{e10}).	
\end{lemma}

Next, we analyze the convergence of (\ref{e10}).

\begin{theorem}\label{t1}
	Under Assumption \ref{a3},   the trajectory $ \left(\boldsymbol{z}(t), \boldsymbol{\lambda}(t), \boldsymbol{\zeta}(t)\right) $ of  (\ref{e10}) is bounded and converges to an equilibrium point of (\ref{e10}),
	namely,  $ \boldsymbol{z}(t) $ converges to a solution of $ \mathrm{VI}(\boldsymbol{\Xi}, \boldsymbol{g}(\boldsymbol{z})) $
	satisfying (\ref{aaf}).
\end{theorem}

\section{Equilibrium Analysis}\label{s6}
In this section,  we   show that an equilibrium obtained from Algorithm 1 induces an $ \epsilon $-GNE of original game (\ref{f5}). Moreover, we describe the  bound related to $ \epsilon $.

The following idea is different from that given
in \cite{chen2021distributed}. As the estimation of $ \epsilon $ is actually reflected by solving GNE of the  approximate problem dependent on dynamics, 	
we consider establishing the relationship between the approximation accuracy and  $ \epsilon $ from the perspective of the nonlinear perturbed system.

Under Assumption \ref{a3},  the pseudo-gradient $ F $ is  strictly monotone with respect to $\boldsymbol{x}$, which implies that $ \boldsymbol{z}^{*} \in \mathrm{SOL}(\boldsymbol{\Xi}, \boldsymbol{g}(\cdot)) $ contains a unique $\boldsymbol{x}^{*}$, but the optimal $\boldsymbol{\sigma}^{*}$ may not be unique. Moreover, if the form of cost function $ J_{i} $ is fixed,
then different polyhedron approximations result in different
variational inequality solutions. Since $ \boldsymbol{\mathcal{P}}_{v} $ determines $ \boldsymbol{\Xi} $, we write
$ \boldsymbol{x}^{*} = \boldsymbol{x}^{*} (\boldsymbol{\mathcal{P}}_{v}) $,  $ \boldsymbol{z}^{*}(\boldsymbol{\mathcal{P}}_{v}) =\operatorname{col}\{ \boldsymbol{x}^{*} (\boldsymbol{\mathcal{P}}_{v}), \boldsymbol{\sigma}^{*} (\boldsymbol{\mathcal{P}}_{v}) \}  $ for game (\ref{f7}). Also, denote $ \boldsymbol{x}^{*} (\boldsymbol{\mathcal{M}}) $ as a GNE of game (\ref{f5}).

Take
\begin{equation}\label{fgg}
	\boldsymbol{\mathcal{P}}_{v_{1}}=\prod_{i=1}^{N} \mathcal{P}_{v_{1, i}}^{i}, \quad \boldsymbol{\mathcal{P}}_{v_{2}}=\prod_{i=1}^{N} \mathcal{P}_{v_{2, i}}^{i}
\end{equation}
as two inscribed polyhedrons of $ \boldsymbol{\mathcal{M}} $.  With the definition of $ \boldsymbol{\mathcal{P}}_{v_{1}} $ and  $ \boldsymbol{\mathcal{P}}_{v_{2}} $, for the $ i $th player,
\begin{equation}\label{bb1}
	\mathcal{P}_{v_{1, i}}^{i}=\left\{\omega_{i} \in \mathbb{R}^{n}: A_{1,i} \omega_{i} \leq d_{1,i}\right\},\;   A_{1,i} \in \mathbb{R}^{q_{1,i}\times n},
\end{equation}
\begin{equation}\label{bb2}
	\mathcal{P}_{v_{2, i}}^{i}=\left\{\omega_{i} \in \mathbb{R}^{n}: A_{2,i} \omega_{i} \leq d_{2,i}\right\},\; A_{2,i} \in \mathbb{R}^{q_{2,i}\times n}.
\end{equation}
Before revealing the $ \epsilon $-relationship of $ \boldsymbol{x}^{*}(\boldsymbol{\mathcal{P}}) $ and $ \boldsymbol{x}^{*}(\boldsymbol{\mathcal{M}}) $
in  game (\ref{f7}) and  original game (\ref{f5}),
we first investigate the  relationship between  $ \boldsymbol{x}^{*}(\boldsymbol{\mathcal{P}}_{v_{1}}) $ and  $ \boldsymbol{x}^{*}(\boldsymbol{\mathcal{P}}_{v_{2}}) $ (i.e., the approximation accuracy between $ \boldsymbol{x}^{*}(\boldsymbol{\mathcal{P}}_{v_{1}}) $ and  $ \boldsymbol{x}^{*}(\boldsymbol{\mathcal{P}}_{v_{2}}) $) of (\ref{f7}).

Define
$\boldsymbol{B}_{1}=\operatorname{Diag}\left(B^{1}_{1}, \ldots, B^{1}_{N}\right) \in \mathbb{R}^{N \times(n N+q_{1})}$, $\boldsymbol{C}_{1}=\operatorname{Diag}\left(C^{1}_{1}, \ldots, C^{1}_{N}\right) \in \mathbb{R}^{nN \times(n N+q_{1})}$, where $ 	B^{1}_{i}=\left[\mathbf{0}_{n}^{\mathrm{T}}, d_{1,i}^{ \mathrm{T}}\right] \in \mathbb{R}^{1 \times\left(n+q_{1,i}\right)} $, $ C_{i}^{1}=\left[-I_{n}, B_{1, i}^{\mathrm{T}}\right] \in \mathbb{R}^{n \times\left(n+q_{1, i}\right)} $ and $ q_{1}=\sum_{i=1}^{N}q_{1,i} $.  $ \boldsymbol{B}_{2} $ and $ \boldsymbol{C}_{2} $ are denoted in a similar way.

Recalling the fact (\ref{fff}) of $\Omega_{i}$, with employing a new variable $\boldsymbol{\xi}=\operatorname{col}\{\xi_{1}, \cdots, \xi_{N}\} \in \mathbb{R}^{n N}$, (\ref{aaf}) on $ \boldsymbol{\mathcal{P}}_{v_{1}} $ is equivalent to
$$
\begin{aligned}
	&	\boldsymbol{0}_{nN}  \in\boldsymbol{g}(\boldsymbol{z}^{*})+\boldsymbol{B}_{1}^{\mathrm{T}} \boldsymbol{\lambda}^{*}+\boldsymbol{C}_{1}^{\mathrm{T}}\boldsymbol{\xi}^{*}+\mathcal{N}_{\boldsymbol{\Phi}}(\boldsymbol{z}^{*}), \\
	&	0 \leq-(\boldsymbol{B}_{1} \boldsymbol{z}^{*}-\boldsymbol{b})^{\mathrm{T}} \cdot \boldsymbol{1}_{N}, \quad 0=(\boldsymbol{B}_{1} \boldsymbol{z}^{*}-\boldsymbol{b})^{\mathrm{T}} \boldsymbol{\lambda}^{*},\\
	&	\boldsymbol{0}_{N}   =L\boldsymbol{\lambda}^{*},\\
	&	\boldsymbol{0}_{nN} = \boldsymbol{C}_{1}\boldsymbol{z}^{*}.
\end{aligned}
$$
%
Let $ \boldsymbol{y} = \operatorname{col}\{\boldsymbol{z},\boldsymbol{\lambda},\boldsymbol{\zeta},\boldsymbol{\xi}\} $, $\boldsymbol{R}=\boldsymbol{\Phi} \times \mathbb{R}_{+}^{N} \times \mathbb{R}^{N} \times \mathbb{R}^{n N}$. Then  Algorithm 1 on  $  \boldsymbol{\mathcal{P}}_{v_{1}} $  is equivalent to
\begin{equation}\label{e13}
	\boldsymbol{y}=D_{\boldsymbol{\mathcal{P}}_{v_{1}}}(\boldsymbol{y}),
\end{equation}
where
$$
D_{\boldsymbol{\mathcal{P}}_{v_{1}}}(\boldsymbol{y})=\!\left[\!\begin{array}{l}
	\Pi_{\boldsymbol{\Phi}}\left(\boldsymbol{z}-\boldsymbol{g}(\boldsymbol{z})-\boldsymbol{B}_{1}^{\mathrm{T}}\boldsymbol{\lambda}-\boldsymbol{C}_{1}^{\mathrm{T}}\boldsymbol{\xi} \right)-\boldsymbol{z}\\
	\left[\boldsymbol{\lambda}+\boldsymbol{B}_{1}\boldsymbol{z}-\boldsymbol{b}-L\boldsymbol{\lambda}-L\boldsymbol{\zeta}\right]^{+}-\boldsymbol{\lambda}\\
	L\boldsymbol{\lambda}\\
	\boldsymbol{C}_{1}\boldsymbol{z}
\end{array}\!\!\right]\!.
$$

From Theorem \ref{t1}, the whole dynamics of system
(\ref{e13}) is globally asymptotically stable. According
to this property,  with  the converse Lyapunov theorem in \cite{khalil2002nonlinear}, there exists a  Lyapunov  function $ V_{\boldsymbol{\mathcal{P}}_{v_{1}}}(\boldsymbol{y}) $  satisfying the following inequalities,
\begin{equation}\label{5q}
	\begin{array}{c}
		\alpha_{1}(\|\boldsymbol{y}-\boldsymbol{y}^{*}(\boldsymbol{\mathcal{P}}_{v_{1}})\|) \leq V_{\boldsymbol{\mathcal{P}}_{v_{1}}}(\boldsymbol{y}) \leq \alpha_{2}(\|\boldsymbol{y}-\boldsymbol{y}^{*}(\boldsymbol{\mathcal{P}}_{v_{1}})\|),\\
		\dot{V}_{\boldsymbol{\mathcal{P}}_{v_{1}}} \leq-\alpha_{3}(\|\boldsymbol{y}-\boldsymbol{y}^{*}(\boldsymbol{\mathcal{P}}_{v_{1}})\|), \\
		\left\|\frac{\partial V_{\boldsymbol{\mathcal{P}}_{v_{1}}}}{\partial \boldsymbol{y}}\right\| \leq \alpha_{4}(\|\boldsymbol{y}-\boldsymbol{y}^{*}(\boldsymbol{\mathcal{P}}_{v_{1}})\|),
	\end{array}
\end{equation}
where $\alpha_{1}$, $\alpha_{2}$, $\alpha_{3}$, $\alpha_{4}$ are  class-$ \mathcal{K} $ functions, $ \boldsymbol{y}^{*}(\boldsymbol{\mathcal{P}}_{v_{1}})=\operatorname{col}\{\boldsymbol{z}^{*}(\boldsymbol{\mathcal{P}}_{v_{1}}),\boldsymbol{\lambda}^{*}(\boldsymbol{\mathcal{P}}_{v_{1}}),\boldsymbol{\zeta}^{*}(\boldsymbol{\mathcal{P}}_{v_{1}}),\boldsymbol{\xi}^{*}(\boldsymbol{\mathcal{P}}_{v_{1}})\} $ is an equilibrium point of (\ref{e13}).

Analogously, the dynamics on  $  \boldsymbol{\mathcal{P}}_{v_{2}} $ is
\begin{equation}\label{5e}
	\boldsymbol{y}=D_{\boldsymbol{\mathcal{P}}_{v_{2}}}(\boldsymbol{y}),
\end{equation}
where
$$
D_{\boldsymbol{\mathcal{P}}_{v_{2}}}(\boldsymbol{y})=\!\left[\!\begin{array}{l}
	\Pi_{\boldsymbol{\Phi}}\left(\boldsymbol{z}-\boldsymbol{g}(\boldsymbol{z})-\boldsymbol{B}_{2}^{\mathrm{T}}\boldsymbol{\lambda}-\boldsymbol{C}_{2}^{\mathrm{T}}\boldsymbol{\xi} \right)-\boldsymbol{z}\\
	\left[\boldsymbol{\lambda}+\boldsymbol{B}_{2}\boldsymbol{z}-\boldsymbol{b}-L\boldsymbol{\lambda}-L\boldsymbol{\zeta}\right]^{+}-\boldsymbol{\lambda}\\
	L\boldsymbol{\lambda}\\
	\boldsymbol{C}_{2}\boldsymbol{z}
\end{array}\!\!\right]\!.
$$

Note that (\ref{5e}) can be regarded as a perturbed system of (\ref{e13}). For clarification, let $ \boldsymbol{\Gamma}= \boldsymbol{z}-\boldsymbol{g}(\boldsymbol{z})$, $ \boldsymbol{\Lambda}=\boldsymbol{\lambda} -\boldsymbol{b}-L\boldsymbol{\lambda}-L\boldsymbol{\zeta}$. Denote
$$ e(\boldsymbol{y})\triangleq D_{\boldsymbol{\mathcal{P}}_{v_{2}}}(\boldsymbol{y})-D_{\boldsymbol{\mathcal{P}}_{v_{1}}}(\boldsymbol{y}),  $$ then (\ref{5e}) is converted into
\begin{equation}\label{aa}
	\boldsymbol{y}=D_{\boldsymbol{\mathcal{P}}_{v_{1}}}(\boldsymbol{y})+e(\boldsymbol{y}),
\end{equation}
where the perturbation term is
$$
e(\boldsymbol{y})\!\!=\!\!\left[\!\begin{array}{c}
	\!\Pi_{\boldsymbol{\Phi}}\left(\boldsymbol{\Gamma}\!-\!\boldsymbol{B}_{2}^{\mathrm{T}}\boldsymbol{\lambda}-\boldsymbol{C}_{2}^{\mathrm{T}}\boldsymbol{\xi} \right)\!\!-\!\Pi_{\boldsymbol{\Phi}}\left(\boldsymbol{\Gamma}\!-\!\boldsymbol{B}_{1}^{\mathrm{T}}\boldsymbol{\lambda}-\boldsymbol{C}_{1}^{\mathrm{T}}\boldsymbol{\xi} \right)\\
	\left[\boldsymbol{\Lambda}+\boldsymbol{B}_{2}\boldsymbol{z}\right]^{+}-\left[\boldsymbol{\Lambda}+\boldsymbol{B}_{1}\boldsymbol{z}\right]^{+}\\
	\boldsymbol{0}\\
	(\boldsymbol{C}_{1}-\boldsymbol{C}_{2})\boldsymbol{z}
\end{array}\!\!\right]\!.
$$
Take $ \boldsymbol{y}^{*}(\boldsymbol{\mathcal{P}}_{v_{2}}) $ as an equilibrium point of (\ref{aa}). After this conversion, we can obtain the upper bound of the approximation accuracy between $ \boldsymbol{y}^{*}(\boldsymbol{\mathcal{P}}_{v_{1}}) $ and  $ \boldsymbol{y}^{*}(\boldsymbol{\mathcal{P}}_{v_{2}}) $ (i.e., $ \boldsymbol{x}^{*}(\boldsymbol{\mathcal{P}}_{v_{1}}) $ and  $ \boldsymbol{x}^{*}(\boldsymbol{\mathcal{P}}_{v_{2}}) $) by investigating  $ e(\boldsymbol{y}) $ between
(\ref{e13}) and (\ref{aa}).

Note that $ e(\boldsymbol{y}) $ reflects the difference in continuous-time projected dynamics on $ \boldsymbol{\mathcal{P}}_{v_{1}} $ and  $ \boldsymbol{\mathcal{P}}_{v_{2}} $, respectively. Recalling the definition of inscribed
polyhedrons in (\ref{bbs}),
$ e(\boldsymbol{y}) $ is  basically affected by different hyperplanes (their corresponding normal vectors and displacement terms) in $ \boldsymbol{\mathcal{P}}_{v_{1}} $ and  $ \boldsymbol{\mathcal{P}}_{v_{2}} $, where the
distance between hyperplanes can be measured by angular
metric.

As defined in (\ref{fgg})-(\ref{bb2}), without losing generality, consider $ q_{1,i} \leq q_{2, i} $. Let  $ A_{2,i}^{l} $ be any row of matrix $ A_{2,i} $, $\forall i \in \mathcal{I}$,  $0\leq l \leq q_{2,i}$, and $ A_{1,i}^{j(l)} $ be the corresponding row of matrix $ A_{1,i} $. Accordingly, denote $\tau^l_i\in[0,\pi/2)$ as the angular metric of $ A_{2,i}^{l} $ and   $ A_{1,i}^{j(l)} $, where
$\tau_{i}^{l}=  \psi\left(A_{2,i}^{l},A_{1,i}^{j(l)}\right)$.
The following lemma gives an upper bound of $ \|e(\boldsymbol{y})\| $.
\begin{lemma}\label{t3}	
	Under Assumption \ref{a3},	on  $ \bar{\boldsymbol{\Omega}}= \boldsymbol{R}\cap\{\|\boldsymbol{y}-\boldsymbol{y}^{*}(\boldsymbol{\mathcal{P}}_{v_{1}})\|<r \} $, the perturbation term
	$ e(\boldsymbol{y}) $  of (\ref{aa}) satisfies
	\begin{equation}\label{e28}
		\|e(\boldsymbol{y})\|\leq\delta=r \sum_{i=1}^{N} q_{i}c_{i}\theta_{i},
	\end{equation}
	where $ c_{i} $  is a finite positive constant,  $ q_{i} =q_{2,i} $ is the number of hyperplanes in 	$ \mathcal{P}_{v_{2, i}}^{i} $, $ \theta_{i}=\max_{0\leq l \leq q_{2,i}} \tau_{i}^{l} $
	for $ i \in \mathcal{I} $.
\end{lemma}
The next lemma explains that  $  \|\boldsymbol{y}^{*}(\boldsymbol{\mathcal{P}}_{v_{1}})- \boldsymbol{y}^{*}(\boldsymbol{\mathcal{P}}_{v_{2}})\| $ is
ultimately bounded by a small bound if $ e(\boldsymbol{y}) $ is small enough,
referring to \cite{khalil2002nonlinear}.

\begin{lemma}\label{fg}
	Take $  V_{\boldsymbol{\mathcal{P}}_{v_{1}}}(\boldsymbol{y}) $  as a Lyapunov function satisfying
	(\ref{5q}) in  set $ \bar{\boldsymbol{\Omega}} $.  Suppose that  	 $ \|e(\boldsymbol{y})\|\leq\delta<\mu \alpha_{3}\left(\alpha_{2}^{-1}\left(\alpha_{1}(r)\right)\right)/\alpha_{4}(r) $,
	with a constant $\mu \in (0,1)$.  Then, for all $ \|\boldsymbol{y}(t_{0})-\boldsymbol{y}^{*}(\boldsymbol{\mathcal{P}}_{v_{1}})\| \leq \alpha_{2}^{-1}(\alpha_{1}(r)) $, the equilibrium $ \boldsymbol{y}^{*}(\boldsymbol{\mathcal{P}}_{v_{2}})  $ of the perturbed system (\ref{aa}) satisfies
	\begin{equation}\label{3q}
		\|\boldsymbol{y}^{*}(\boldsymbol{\mathcal{P}}_{v_{1}})-\boldsymbol{y}^{*}(\boldsymbol{\mathcal{P}}_{v_{2}})\| \leq \rho(\delta)=\alpha_{1}^{-1}\!\left(\!\alpha_{2}\left(\!\alpha_{3}^{-1}\left(\!\frac{\delta \alpha_{4}(r)}{\mu}\right)\right)\right).
	\end{equation}
\end{lemma}

Due to the analysis in Lemma \ref{fg},
for any arbitrarily
small perturbations, there always exists a finite $r$ to satisfy (\ref{3q}).
Clearly, a lower metric yields a lower bound.  It can be
regarded as the robustness of the nominal system with a stable
equilibrium, since arbitrarily small perturbations will not cause
a significant deviation.  Moreover, it follows from (\ref{3q}) that $$ \|\boldsymbol{x}^{*}(\boldsymbol{\mathcal{P}}_{v_{1}})-\boldsymbol{x}^{*}(\boldsymbol{\mathcal{P}}_{v_{2}})\| \leq \rho(\delta). $$
Since $ \rho(0)=0 $ and $  \rho $ is strictly increasing in $ [0, \infty) $, $ \|\boldsymbol{x}^{*}(\boldsymbol{\mathcal{P}}_{v_{1}})-\boldsymbol{x}^{*}(\boldsymbol{\mathcal{P}}_{v_{2}})\|  $ tends to zero as $\delta $ vanishes.
\begin{remark}
	Compared with the analysis in \cite{chen2021distributed}, Lemma \ref{t3}
	does not rely on the Hausdorff metric, which leads to technical
	difficulties in estimating the parameter changes of different
	polyhedrons, and thus can not describe the relationship between
	the approximate accuracy of different polyhedrons and
	the difference between the corresponding equilibria. Instead, by
	introducing angular metric, these difficulties are solved,
	and the upper
	bound of the difference between equilibria can be obtained, which extends the result in \cite{chen2021distributed} and ensures the estimation of
	$ \epsilon $ in the sequel.
\end{remark}

With Lemma \ref{fg}, we finally show that an equilibrium $ \boldsymbol{x}^{*}(\boldsymbol{\mathcal{P}}_{v}) $
obtained from Algorithm 1 induces an $ \epsilon $-GNE of original game
(\ref{f5}) and estimate the approximation accuracy of $ \epsilon $.
\begin{theorem}\label{t9}	
	Under Assumption \ref{a3},
	
	(i) the variational  GNE $ \boldsymbol{x}^{*}(\boldsymbol{\mathcal{M}}) $  of game (\ref{f5}) in the worst case exists and is unique;
	
	(ii) $ \boldsymbol{x}^{*}(\boldsymbol{\mathcal{P}}_{v}) $ of the equilibrium in Algorithm 1 induces an
	$ \epsilon $-GNE of game (\ref{f5}) in the worst case;
	
	(iii) the value of  $\epsilon$  satisfies
	\begin{equation}\label{e25}	
		\begin{aligned}
			\epsilon&\leq2\varsigma_{i} \alpha_{1}^{-1}\left(\alpha_{2}\left(\alpha_{3}^{-1}\left(\frac{\delta \alpha_{4}(r)}{\mu}\right)\right)\right),
		\end{aligned}
	\end{equation}
	where   the  constant $ \mu \in(0,1)  $, $\alpha_{1}$, $\alpha_{2}$, $\alpha_{3}$, $\alpha_{4}$ are  class-$ \mathcal{K} $ functions in (\ref{5q}), $ \varsigma_{i} $ is the Lipschitz constant of $ J_{i} $. Specifically,
	\begin{align}\label{e36}
		&\delta
		=r \sum_{i=1}^{N} \frac{q_{i}c_{i}}{\sqrt{\frac{2}{h_{i}\nu_{i}}-1}},
	\end{align}
	where  $h_{i}=H(\mathcal{P}^{i}_{v_{i}},  \mathcal{M}_{i})$  is  the Hausdorff  distance between $ \mathcal{P}^{i}_{v_{i}}$ and $ \mathcal{M}_{i}$, $q_i$ is the number of hyperplanes in $ \mathcal{P}^{i}_{v_{i}}$, $ \nu_{i} $ is a
	constructive curvature related merely to the structure of $ \mathcal{M}_{i} $,
	$ c_i $ is a finite constant.
\end{theorem}

From Theorem \ref{t9}, the upper bound of $ \epsilon $ is proportional
to the bound of $ \delta $. With the expression of $ \delta $ in (\ref{e36}), when constructing
polyhedrons with more vertices, we obtain more hyperplanes
enclosed the polyhedrons (more rows of matrix $ A_i $ and vectors
$ d_i $), which results in a lower metric and higher accuracy of $ \epsilon $. Actually,
there are developed investigations on how to construct
a proper inscribed polyhedron \cite{kamenev1992class,kamenev2003self}.  When the vertices or faces are constructed successively, we can find a proper inscribed polyhedron by the iterative algorithms based on Hausdorff metric. The main idea of iterative algorithms is to construct a polyhedron $\mathcal{P}_{v(k+1)}=\operatorname{conv}\left(\mathcal{P}_{v(k)} \cup\left\{w_{k+1}\right\}\right)$
	every iteration, where $ v(k) $ is the number of vertices in $ \mathcal{P}_{v(k)} $, $ w_{k+1} $ is a point from  $\partial \mathcal{M}$
	(i.e., the boundary of $ \mathcal{M} $). The Hausdorff metric satisfies $ H(\mathcal{P}_{v(k)}, \mathcal{M}) \leq C_{\mathcal{M}}\cdot{v(k)^{(n-1)/2}} $,
	where $ C_{\mathcal{M}} $ is a constant related with the curvature of $ \mathcal{M} $.
	One of the methods of constructing point $ w_{k+1} $
	is described as follows.	For $ u\in \mathbb{R}^{n} $, denote $g_{\mathcal{M}}(u)=\max \{\langle u, x\rangle: x \in \mathcal{M}\}$ as the support function of $ \mathcal{M} $ on the unit sphere of directions $S^{n-1}=\left\{u \in \mathbb{R}^{n}:\|u\|=1\right\}$. 	The additional point $  w_{k+1} \in\partial \mathcal{M} $
	belongs to the support plane parallel to the hyperplane in $  \mathcal{P}_{v_{k}} $, for which the quantity $ g_{\mathcal{M}}(u)- g_{\mathcal{P}_{v(k)}}(u)$  attains its maximum on the set of external normals $  u \in S^{n-1} $ to the hyperplanes of $ \mathcal{P}_{v(k)}$. The initial polyhedron could be constructed by  the method \cite{2008Approximation}.
In addition, since the parameter set constraint of each
player is private information to itself, different players can
approximate their parameter sets through different construction
methods separately, in advance and offline.

\section{Numerical experiments}\label{s7}
In this section, we examine the  approximation accuracy of Algorithm 1 on demand response management problems  under uncertainty as in \cite{ye2016game,wei2014energy}.

Consider a game with $ N = 10 $ electricity users with the demand of energy consumption. For $ i \in \mathcal{I} =
{1, \cdots ,10} $,  $ x_{i}\in \Theta_i$ is the energy consumption of the $ i $th user, where
$ \Theta_i=\{ x_i\in \mathbb{R}^{2} : c_{1}\boldsymbol{1}_{2}\leq x_i \leq c_{2}\boldsymbol{1}_{2} \} $ with $ c_{1}=-15 $, $ c_{2}=20 $.
In this network  game, each user needs
to solve the following  problem given the other
users’ profile $ \boldsymbol{x}_{-i} $,
\begin{equation}
	\begin{aligned}
		&\min \limits_{x_{i} \in \Theta_i}  \frac{1}{2}(x_{i}-\varpi_{i})^{\mathrm{T}}(x_{i}-\varpi_{i})
		-x_{i}^{\mathrm{T}}p(\mathcal{Q}(\boldsymbol{x})), \\
		\text { s.t. } & \sum_{j=1}^{N} a^{\mathrm{T}} x_{j} \!\leq\! b, \quad a\! \in\! \mathbf{E}_{(3,2)}(2,2),\quad \forall j\in \mathcal{I},
	\end{aligned}
\end{equation}
where $ \varpi_{i}=(5-i)\boldsymbol{1}_{2}\in \mathbb{R}^{2} $ is the nominal value of energy consumption,
and   $p=N(\boldsymbol{1}_{2}-\mathcal{Q}(\boldsymbol{x})) $ is the pricing function with $
\mathcal{Q}(\boldsymbol{x})=\frac{1}{N} \sum_{j=1}^{N} x_{j}
$ as an aggregative term.  All  electricity users  need to meet the  coupled
inequality constraint with   the parameter  $ a \in \mathbb{R}^{2} $ satisfying an elliptical
region
$$
\mathbf{E}_{3,2}(2,2)=\left\{a \in \mathbb{R}^{2}: \frac{(a_{ 1}-2)^{2}}{3^{2}}+\frac{(a_{ 2}-2)^{2}}{2^{2}} \leq 1\right\} .
$$
Take a ring graph as the communication network $ \mathcal{G} $,
$$
1\rightleftarrows2\rightleftarrows\cdots\rightleftarrows10\rightleftarrows1.
$$
Meanwhile, we set tolerance as $ t_{tol} = 10^{-4} $
and the terminal criterion as $\|\dot{\boldsymbol{y}}(t)\| \leq t_{tol }$.
We employ inscribed rectangles to approximate $ \mathbf{E}_{(3,2)}(2,2) $, where  the trajectories of one dimension of each $ x_i $ are shown in Fig. 1. Then we verify the approximation accuracy of Algorithm 1. We approximate $ \mathbf{E}_{(3,2)}(2,2) $ with inscribed triangles, rectangles, hexagons,
octagons, decagons, and dodecagons, respectively. Fig. 2 presents different strategy trajectories  of one fixed player with different approximations.  The vertical axis represents the value of the convergent $ \epsilon $-GNE and the horizontal axis
represents the real running time of Algorithm 1. The results imply that when we choose a more accurate approximation, equilibria with different polyhedrons
get closer to the exact solution.
\begin{figure}
	\centering	
	\includegraphics[width=0.6\linewidth]{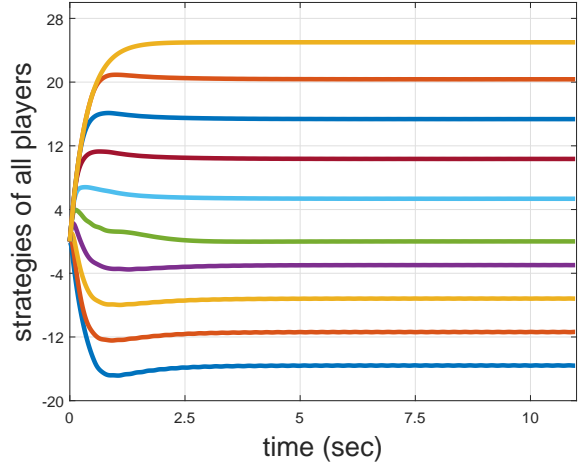}\\
	\caption{ Trajectories of all players’ strategies.}
	\label{fig4}
\end{figure}

\begin{figure}
	\centering	
	\includegraphics[width=0.6\linewidth]{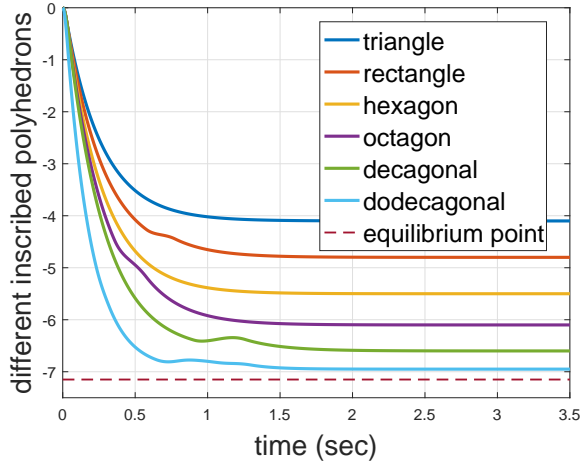}\\
	\caption{ Trajectories of approximation by different inscribed polyhedrons.}
	\label{fig6}
\end{figure}

Additionally, recalling the definition of $ \epsilon $-GNE, the numerical
values of $ \epsilon $ under different types of approximation are listed in Table 1.
Obviously, the value of $ \epsilon $ decreases
with the increase of the vertices of polyhedrons and the decrease of Hausdorff distances, which is consistent
with the approximation results.

\begin{figure}
	\centering	
	\includegraphics[width=0.6\linewidth]{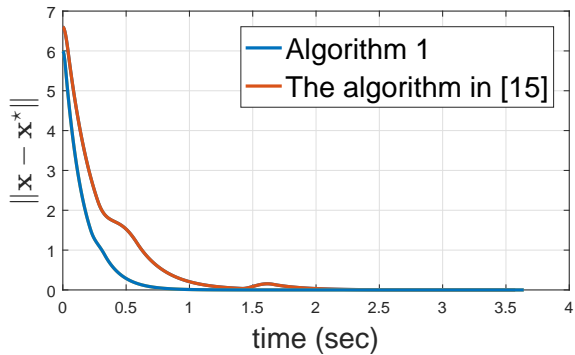}\\
	\caption{ The comparison of the performance of our algorithm and the
		algorithm in \cite{chen2021distributed}.}
	\label{fig7}
\end{figure}

We further verify the effectiveness of our algorithm by comparing it with the algorithm of \cite{chen2021distributed}. Fig. \ref{fig7} shows comparative results for our algorithm and the method proposed in \cite{chen2021distributed}. The results imply that both of them are convergent, and  (\ref{e10}) is with a faster convergence rate because (\ref{e10}) has lower dimensions and less complexity.
\begin{table*}
	\centering
	\label{tab1}
	\caption{Performance of different approximations.}
	\setlength\tabcolsep{7.58pt}
	\renewcommand\arraystretch{1}
	\begin{tabular}{l|c|c|c|c|c|c}
		\hline
		\hline
		\specialrule{0em}{0.2pt}{0.5pt}
		Polyhedrons &  Triangle &  Rectangle &  Hexagon &  Octagon&Decagonal&Dodecagonal \\ \hline
		Values of $ \epsilon $ & 16.0416 & 11.8262 & 6.6113 & 3.9556 &1.5406&0.7054\\ \hline\hline
	\end{tabular}
\end{table*}

\section{Conclusion}\label{s8}
A distributed game with coupled inequality constraints has been studied in this paper, where parameters in constraints are from general uncertain convex sets. By employing inscribed polyhedrons
to approximate parameter sets,  a distributed algorithm has been proposed  for seeking an $ \epsilon $-GNE in the worst case, and the convergence of the algorithm has been shown. With the help of convex set geometry and metric spaces, the approximation accuracy affected by  different inscribed polyhedrons is analyzed. Moreover, with the proof that the equilibrium point of the algorithm is an $ \epsilon $-GNE of the original problem, an upper bound
of the value of $ \epsilon $ has been estimated by analyzing a perturbed system.


%
%
%
%
%
%
%
%
\renewcommand\bibname{References}
\bibliographystyle{IEEEtran}
\bibliography{references}

\begin{appendices}
\section{Proof of Lemma  \ref{l3}}\label{a13}
	(i) Consider $(\boldsymbol{z}^{*},\boldsymbol{\lambda}^{*},\boldsymbol{\zeta}^{*})$ as an equilibrium point of (\ref{e10}).
By properties of  normal cones to  nonempty closed convex sets,  at the equilibrium point, $ \dot{\boldsymbol{z}}= \boldsymbol{0}_{nN} $ implies that $\Pi_{\boldsymbol{\Omega}}(\boldsymbol{z}^{*}-$ $\boldsymbol{g}(\boldsymbol{z}^{*})-\boldsymbol{B}^{\mathrm{T}}\boldsymbol{\lambda}^{*})=\boldsymbol{z}^{*}$. Then it follows from Lemma 2.38 of \cite{ruszczynski2011nonlinear} that  $ -\boldsymbol{g}(\boldsymbol{z}^{*})-\boldsymbol{B}^{\mathrm{T}} \boldsymbol{\lambda}^{*}\in \mathcal{N}_{\boldsymbol{\Omega}}(\boldsymbol{z}^{*}).$

Moreover, we set $ \dot{\boldsymbol{\zeta}}=\boldsymbol{0}_{N} $ and $ \dot{\boldsymbol{\lambda}}=\boldsymbol{0}_{N} $, which obtain $ L\boldsymbol{\lambda}^{*}=\boldsymbol{0}_{N} $ and
$
\boldsymbol{B}\boldsymbol{z}^{*}-\boldsymbol{b}-L\boldsymbol{\zeta}^{*}\in \mathcal{N}_{\mathbb{R}^{N}_{+}}(\boldsymbol{\lambda}^{*}).
$
It implies that $
\boldsymbol{B}\boldsymbol{z}^{*}-\boldsymbol{b}-L\boldsymbol{\zeta}^{*}\leq \boldsymbol{0}_{N}.
$
Because the graph $ \mathcal{G} $ is undirected and connected, $ \boldsymbol{1}_{N}^{\mathrm{T}}L= \boldsymbol{0}_{N}^{\mathrm{T}}$, and $
\boldsymbol{1}_{N}^{\mathrm{T}}(\boldsymbol{B}\boldsymbol{z}^{*}-\boldsymbol{b})\leq 0.
$
Also, take $ \boldsymbol{\varsigma}^{*} \in \mathcal{N}_{\mathbb{R}_{+}^{N}}(\boldsymbol{\lambda}^{*}) $. Then we have $
\boldsymbol{B}\boldsymbol{z}^{*}-\boldsymbol{b}-L\boldsymbol{\zeta}^{*}-\boldsymbol{\varsigma}^{*}=\boldsymbol{0}_{N}
$. When $ \boldsymbol{\lambda}^{*}>\boldsymbol{0}_{N} $, $ \boldsymbol{\varsigma}^{*} = \boldsymbol{0}_{N} $. Then it derives that $ (\boldsymbol{B} \boldsymbol{z}^{*}-\boldsymbol{b})^{\mathrm{T}} \boldsymbol{\lambda}^{*}=\boldsymbol{0}_{N}  $.
When $ \boldsymbol{\lambda}^{*}=\boldsymbol{0}_{N} $,  $ \boldsymbol{\varsigma}^{*} \in -\mathbb{R}^{N}_{+}  $, and $ (\boldsymbol{B} \boldsymbol{z}^{*}-\boldsymbol{b})^{\mathrm{T}} \boldsymbol{\lambda}^{*}=\boldsymbol{0}_{N}  $ is still hold.
Thus, $ \boldsymbol{z}^{*} $ is a variational GNE of game (\ref{f7}).

(ii) When $ \boldsymbol{z}^{*} $ is a variational GNE of game (\ref{f7}), there
exists  $ \boldsymbol{\lambda}^{*} \in \mathbb{R}_{+} $ such that the first order conditions (\ref{aaf}) are
satisfied. It is clear that $ -\boldsymbol{g}(\boldsymbol{z}^{*})-\boldsymbol{B}^{\mathrm{T}} \boldsymbol{\lambda}^{*}\in \mathcal{N}_{\boldsymbol{\Omega}}(\boldsymbol{z}^{*})$ is
equivalent to $\Pi_{\boldsymbol{\Omega}}\left(\boldsymbol{z}^{*}-\boldsymbol{g}(\boldsymbol{z}^{*})-\boldsymbol{B}^{\mathrm{T}}\boldsymbol{\lambda}^{*} \right)=\boldsymbol{z}^{*}$. Furthermore, since $ 0\geq
(\boldsymbol{B}\boldsymbol{z}^{*}-\boldsymbol{b})^{\mathrm{T}}\cdot\boldsymbol{1}_{N}
$, there exists an
$ \boldsymbol{\gamma}\in \mathbb{R}_{+}^{N} $  such that $ 0= (\boldsymbol{B}\boldsymbol{z}^{*}-\boldsymbol{b}+\boldsymbol{\gamma})^{\mathrm{T}}\cdot\boldsymbol{1}_{N}$. Note that $ L\boldsymbol{1}_{N} = \boldsymbol{0}_{N} $ implies $\operatorname{ker}(L)=\operatorname{span}\left\{\mathbf{1}_{N}\right\}$. With $\mathbb{R}^{N}=\operatorname{ker}(L) \oplus \operatorname{Im}(L)$, there
exists $\boldsymbol{\zeta}^{*} \in \operatorname{Im}(L)$ such that $L \boldsymbol{\zeta}^{*}=\boldsymbol{B} \boldsymbol{z}^{*}-\boldsymbol{b}+\boldsymbol{\gamma}$, which
implies $\boldsymbol{B} \boldsymbol{z}^{*}-\boldsymbol{b}-L \boldsymbol{\zeta}^{*} \in \mathcal{N}_{\mathbb{R}_{+}^{N}}\left(\boldsymbol{\lambda}^{*}\right)$. Therefore, $(\boldsymbol{z}^{*},\boldsymbol{\lambda}^{*},\boldsymbol{\zeta}^{*})$
is an equilibrium point of (\ref{e10}).

\section{Proof of Theorem \ref{t1}}\label{a11}
Let $\widehat{\boldsymbol{\Omega}} \triangleq \boldsymbol{\Omega} \times \mathbb{R}_{+}^{N} \times \mathbb{R}^{N}$ and $ s=\operatorname{col}\{\boldsymbol{z},\boldsymbol{\lambda},\boldsymbol{\zeta}\} $.	Define
$$
\widehat{F}(s) \triangleq  \left(\begin{array}{c}
	\boldsymbol{g}(\boldsymbol{z})+\boldsymbol{B}^{\mathrm{T}}\boldsymbol{\lambda} \\
	-\boldsymbol{B}\boldsymbol{z}+\boldsymbol{b}+L\boldsymbol{\lambda}+L\boldsymbol{\zeta} \\
	-L\boldsymbol{\lambda}
\end{array}\right) , $$
$$ U(s)\triangleq\Pi_{\widehat{\Omega}}(s-\widehat{F}(s)).$$

Take the following Lyapunov function
\begin{equation}\label{ss}
	V(t)=-\langle \widehat{F}(s), U(s)-s\rangle-\frac{1}{2}\|U(s)-s\|_{2}^{2}+\frac{1}{2}\left\|s-s^{*}\right\|_{2}^{2},
\end{equation}	
where $ s^{*}=\operatorname{col}\{\boldsymbol{z}^{*},\boldsymbol{\lambda}^{*},\boldsymbol{\zeta}^{*}\} $.
It follows from \cite{fukushima1992equivalent} that $ -\langle \widehat{F}(s), U(s)-s\rangle-\frac{1}{2}\|U(s)-s\|_{2}^{2}\geq 0 $. Thus,
$ V(t)\geq \frac{1}{2}\left\|s-s^{*}\right\|_{2}^{2}\geq 0$, and $  V(t)=0 $ if and only if $ s=s^{*} $.
Moreover, referring to \cite{liang2017distributed}, $ \dot{V}(t) $ can be calculated as
\begin{equation}\label{www}
	\begin{aligned}
		\dot{V}(t)
		&\leq -(\widehat{F}(s)-\widehat{F}(s^{*}))^{\mathrm{T}}(s-s^{*})\\
		&= -\left(\boldsymbol{z}-\boldsymbol{z}^{*}\right)^{\mathrm{T}}\left(g(\boldsymbol{z})-g(\boldsymbol{z}^{*}) \right)-\boldsymbol{\lambda}^{\mathrm{T}}L \boldsymbol{\lambda}.
	\end{aligned}
\end{equation}	
Due to the monotonicity of  $ g(\boldsymbol{z}) $,
it derives that $ \dot{V}(t)\leq 0 $. Hence, the trajectory of algorithm (\ref{e10}) is bounded and  any finite equilibrium
point  of (\ref{e10}) is Lyapunov stable.

Furthermore, denote the set of points satisfying $ \dot{V}(t)= 0 $ by
$
E_{v} \triangleq\left\{(\boldsymbol{z}, \boldsymbol{\lambda},\boldsymbol{\zeta}): \dot{V}(t)=0)\right\}.
$
From (\ref{www}), there holds
\begin{equation}\label{la}
	E_{v} \subseteq\left\{(\boldsymbol{z}, \boldsymbol{\lambda},\boldsymbol{\zeta}): \boldsymbol{z}=\boldsymbol{z}^{*}, L\boldsymbol{\lambda}=0 \right\}.    	
\end{equation}
Then we claim that the maximal invariance set $ R $ within the
set $ E_{v} $ is exactly the equilibrium point of (\ref{e10}). It follows from the invariance principle (Theorem 4.4 of \cite{khalil2002nonlinear}) that $ (\boldsymbol{z}(t), \boldsymbol{\lambda}(t),\boldsymbol{\zeta}(t))\rightarrow R$ as $ t\rightarrow\infty $, and $ R $ is a positive invariant set. Consider a trajectory $ (\bar{\boldsymbol{z}},  \bar{\boldsymbol{\lambda}}, \bar{\boldsymbol{\zeta}}) $ in $ R $.
Note that (\ref{la}) implies $ \dot{\bar{\boldsymbol{z}}}=\boldsymbol{0} $, $ \dot{\bar{\boldsymbol{\zeta}}}=\boldsymbol{0} $, and $ \dot{\bar{\boldsymbol{\lambda}}}=constant $.   Due to the boundness of the trajectory,  it leads   to a contradiction if $ \dot{\bar{\boldsymbol{\lambda}}}\neq \boldsymbol{0} $. Hence, any point in $ R $ is an equilibrium point of algorithm (\ref{e10}). By Corollary 4.1 in \cite{khalil2002nonlinear},  system (\ref{e10}) converges to its equilibrium point. Therefore, based on Lemma \ref{l3},  $ \boldsymbol{z}(t) $ converges to a solution of $ \mathrm{VI}(\boldsymbol{\Xi}, \boldsymbol{g}(\boldsymbol{z})) $
satisfying (\ref{aaf}).

\section{Proof of Lemma \ref{t3}}\label{a16}

We will prove the conclusion of Lemma \ref{t3} in two steps.

Step 1: Denote  $\mathcal{P}_{v_{1}}^{1}=\left\{\omega \in \mathbb{R}^{n}: A_{1} \omega \leq d_{1}\right\}$ as  an inscribed polyhedron of  a convex and compact set $ \mathcal{M} \subseteq \mathbb{R}^{n} $ with $ V_{1} $ as  the set of
vertices on the boundary of $ \mathcal{M} $. Take  $ \mathcal{P}_{v_{2}}^{2}=\left\{\omega \in \mathbb{R}^{n}: A_{2} \omega \leq d_{2}\right\} $ as another inscribed polyhedron whose vertices consist of $ V_{2}= V_{1}\cup \{v_{0}\} $, with $ v_{0} $ as an additional vertex on the boundary of $ \mathcal{M} $.
We first prove that 	for $ A_{2}^{l} $ as any row of matrix $ A_{2} $,  there exists   a  corresponding row  $ A_{1}^{j(l)} $ of matrix $ A_{1} $ such that  	
$ 	\left\|A^{l}_{2}-A^{j(l)}_{1}\right\| \leq c\tau^{l} $,
where $ \tau^{l}=\psi\left(A_{2}^{l},A_{1}^{j(l)}\right) \in  [0,\pi/2)$ is the angular metric  between $ A_{2}^{l} $ and $ A^{j(l)}_{1} $, $ c $ is a finite positive constant.

Suppose that there are $ q_{1} $ rows of $ A_{1} $ and $ d_{1} $, $ q_{2} $ rows of $ A_{2} $ and $ d_{2} $, the first $ q_{1}-1 $ rows of $ A_{1} $ are the same as the first $ q_{1}-1 $ rows of $ A_{2} $.  Thus, we only need to investigate the difference between $ A_{1}^{q_{1}} $ and the last $ q_{2}-q_{1}+1 $ rows of $ A_{2} $.

Note that the dimension of each   hyperplane is $ n-1 $, and  normalized vectors $ A_{2}^{l} $ (or $ A_{1}^{j(l)} $) represent normal vectors of hyperplanes enclosing the polyhedron $ \mathcal{P}_{v_{2}}^{2} $ (or $ \mathcal{P}_{v_{1}}^{1} $). It follows from Lemma  \ref{l5} that the angle between
two hyperplanes uniquely equals to that between their normal
vectors. Then there exists a derived angular metric   and a
corresponding scalar $ \tau^{l} \in  [0,\pi/2)$ for $ q_{1} \leq l \leq q_{2}$ such that
$\tau^{l}  =\psi\left(A_{2}^{l},A_{1}^{q_{1}}\right) $.

Additionally, referring to \cite[Theorem 2.21]{kato2013perturbation}, there exists a
derived gap metric $ \upsilon\left(A_{2}^{l},A_{1}^{q_{1}}\right) $
such that
$$
\left\|A_{2}^{l}-A_{1}^{q_{1}}\right\| \leq \frac{1+\left\|A_{1}^{q_{1}}\right\|^{2}}{\sqrt{1+\left\|A_{1}^{q_{1}}\right\|^{2}}-1} \cdot v\left(A_{2}^{l}, A_{1}^{q_{1}}\right).
$$
According to the definition of the gap metric in \cite{qiu2005unitarily} and \cite{el1985gap},
there holds
$$
v\left(A_{2}^{l}, A_{1}^{q_{1}}\right)=\sin \tau^{l}.
$$
Since $A_{1} $ and  $A_{2} $ are with normalized rows, 	with the fact
$ 	\sin \tau^{l}\leq \tau^{l}   $, there exists a constant $ c $ such that
$
\left\|A_{2}^{l}-A_{1}^{q_{1}}\right\|
\leq c\tau^{l} .\\
$	

Step 2:   Take $ \boldsymbol{\mathcal{P}}_{v_{1}} $ and $ \boldsymbol{\mathcal{P}}_{v_{2}} $ defined in (\ref{bb1}) and (\ref{bb2}) as  two arbitrarily inscribed polyhedrons of $ \boldsymbol{\mathcal{M}} $. Without losing generality,
consider $ q_{1, i} \leq q_{2, i} $, $ \forall i\in \mathcal{I} $.
If $ q_{1, i} < q_{2, i} $, then we increase the number of the hyperplane in $ \mathcal{P}_{v_{1, i}}^{i} $ successively. The newly added hyperplanes are the same as the $ q_{1, i} $-th hyperplane. Continue this process until $ q_{1, i} = q_{2, i} $.

According to the Lipschitz continuous of the projection,
$$ 	\|e(\boldsymbol{y})\| \leq \|\boldsymbol{B}_{1}-\boldsymbol{B}_{2}\|(\|\boldsymbol{\lambda}\|+\|\boldsymbol{z}\|)+\|\boldsymbol{C}_{1}-\boldsymbol{C}_{2}\|(\|\boldsymbol{\xi}\|+\|\boldsymbol{z}\|)\\. $$
For $ i \in \mathcal{I} $, since $ 	B^{1}_{i}-B^{2}_{i}=\left[\mathbf{0}_{n}^{\mathrm{T}}, \left(d_{1,i }-d_{2,i }\right)^{\mathrm{T}}\right] \in \mathbb{R}^{1 \times\left(n+q_{2,i}\right)}  $ and $ 	C^{1}_{i}-C^{2}_{i}=\left[\!-\boldsymbol{0}_{n\times n}, \left(A_{1,i }\!-\!A_{2,i }\right)\!^{\mathrm{T}}\!\right] $ $\in \mathbb{R}^{n \times\left(n+q_{2,i}\right)} $,
we only need to investigate  $ \|A_{1,i }-A_{2,i }\| $ and $ \|d_{1,i }-d_{2,i }\| $.

It follows from Step 1  that
$	
\left\|A_{2,i}^{l}-A_{1,i}^{j(l)}\right\| \leq \tau_{i}^{l}c_{A,i} \leq \theta_{i}c_{A,i}$, $ \forall i\in \mathcal{I} $,
where $ c_{A,i} $ is a constant for  $ i\in \mathcal{I} $.
Then, $ \|A_{1,i }-A_{2,i }\|\leq q_{i}c_{A,i}\theta_{i}, $
where $ q_{i}=q_{2,i} $ is the number of hyperplanes in 	$ \mathcal{P}_{v_{2, i}}^{i} $. Correspondingly, $ \left\|d_{2,i}-d_{1,i}\right\| \leq q_{i}c_{d,i}\theta_{i}$, where $ c_{d,i} $ is a constant for  $ i\in \mathcal{I} $.  The analysis of other players
is similar to that of player $ i $. To sum up, there exists a finite
constant $ c $ such that $ \|e(\boldsymbol{y})\| $ can be bounded by $\delta$
%
on  $ \bar{\boldsymbol{\Omega}} $, that is,
$$
\small
\begin{aligned}
	\|e(\boldsymbol{y})\|
	&\leq r(\|\boldsymbol{B}_{1}-\boldsymbol{B}_{2}\|+\|\boldsymbol{C}_{1}-\boldsymbol{C}_{2}\|)\\
	& \leq r\sum_{i=1}^{N}\|A_{1,i }-A_{2,i }\|+\|d_{1,i }-d_{2,i }\|\\
	&= r\sum_{i=1}^{N}q_{i}(c_{A,i}+c_{d,i})\theta_{i}=r\sum_{i=1}^{N}q_{i}c_{i}\theta_{i}.	\\
\end{aligned}
$$

\section{Proof of Theorem \ref{t9}}\label{a14}
We first verify the existence and uniqueness of $ \boldsymbol{x}^{*}(\boldsymbol{\mathcal{M}}) $.

Referring to \cite{2008Approximation}, 	for a  convex set $\mathcal{M}$, there exists an  inscribed polyhedron $ \mathcal{P}_{v} $ of $\mathcal{M}$ such that the upper bound of the Hausdorff metric between $\mathcal{M}$ and $ \mathcal{P}_{v} $ satisfies $ H(\mathcal{P}_{v}, \mathcal{M}) \leq C_{\mathcal{M}}\cdot{v^{(n-1)/2}} $,
where $ C_{\mathcal{M}} $ is a constant related with the curvature of $ \mathcal{M} $, and $ v $ is the number of vertices in $ \mathcal{P}_{v} $.  That is to say,  $\lim\limits_{v\rightarrow \infty}H(\boldsymbol{\mathcal{P}}_{v}, \boldsymbol{\mathcal{M}})=0 $. Meanwhile, following from \cite[Lemma 4]{xu2021efficient}, there holds
\begin{equation}\label{aad}
	\tau^{l}_{i}\leq \theta_{i}\leq \frac{1}{\sqrt{\frac{2}{h_{i}\nu_{i}}-1}},
\end{equation}
where  $h_{i}=H(\mathcal{P}^{i}_{v_{1,i}},\mathcal{P}^{i}_{v_{2,i}}) $ represents the Hausdorff  distance between $ \mathcal{P}^{i}_{v_{1,i}}$ and $ \mathcal{P}^{i}_{v_{2,i}}$, $ \nu_{i} $ is a constructive curvature related merely to the structure of
$ \mathcal{M}_{i}$ for $ i\in \mathcal{I}$. Denote $ \boldsymbol{H}=\operatorname{col}\{h_1,\cdots,h_N\} $.  By substituting (\ref{aad}) into (\ref{e28}) and (\ref{3q}), $ \|\boldsymbol{x}^{*}(\boldsymbol{\mathcal{P}}_{v_{1}})-\boldsymbol{x}^{*}(\boldsymbol{\mathcal{P}}_{v_{2}})\| \rightarrow 0 $ as $\boldsymbol{H}\rightarrow0  $,  which means that $ \boldsymbol{x}^{*}(\boldsymbol{\mathcal{P}}_{v}) $ is continuous in $ \boldsymbol{\mathcal{P}}_{v} $ under Hausdorff metric.
Therefore,  there exist a unique $ \boldsymbol{x}^{*}(\boldsymbol{\mathcal{M}}) $ such that
\begin{equation}\label{dee}
	\lim\limits_{v\rightarrow \infty}\boldsymbol{x}^{*}(\boldsymbol{\mathcal{P}}_{v})=\boldsymbol{x}^{*}(\boldsymbol{\mathcal{M}}).
\end{equation}	

Next, we prove that  $  \boldsymbol{x}^{*}(\boldsymbol{\mathcal{P}}_{v}) $ of  approximate
game (\ref{f7}) is an $ \epsilon $-GNE of the original game (\ref{f5}) and estimate $ \epsilon $. Rewrite  $ \delta $  as $ \delta(\boldsymbol{\mathcal{P}}_{v_{1}},\boldsymbol{\mathcal{P}}_{v_{2}}) $.
When  $ \boldsymbol{\mathcal{P}}_{v_{2}} $ is fixed,  $ \delta(\boldsymbol{\mathcal{P}}_{v_{1}},\boldsymbol{\mathcal{P}}_{v_{2}}) $ is continuous in $ \boldsymbol{\mathcal{P}}_{v_{1}} $.
By substituting $ \boldsymbol{\mathcal{P}}_{v_{1}} $ with $ \boldsymbol{\mathcal{P}}_{v_{k}} $, we have
$$ \| \lim\limits_{k\rightarrow \infty}\boldsymbol{x}^{*}(\boldsymbol{\mathcal{P}}_{v_{k}})-\boldsymbol{x}^{*}(\boldsymbol{\mathcal{P}}_{v_{2}})\| \leq \rho(\delta(\lim\limits_{k\rightarrow \infty}\boldsymbol{\mathcal{P}}_{v_{k}},\boldsymbol{\mathcal{P}}_{v_{2}})). $$
Note that (\ref{dee}) is equivalent to
$\lim\limits_{k\rightarrow \infty}\boldsymbol{x}^{*}(\boldsymbol{\mathcal{P}}_{v_{k}})=\boldsymbol{x}^{*}(\boldsymbol{\mathcal{M}}). $
With  $ \lim\limits_{k\rightarrow \infty}\boldsymbol{\mathcal{P}}_{v_{k}}=\boldsymbol{\mathcal{M}}, $
we have
$$
\left\|\boldsymbol{x}^{*}\left(\boldsymbol{\mathcal{P}}_{v}\right)-\boldsymbol{x}^{*}(\boldsymbol{\mathcal{M}})\right\|\leq  \rho(\delta(\boldsymbol{\mathcal{P}}_{v},\boldsymbol{\mathcal{M}})).
$$
Moreover, since $ \delta(\boldsymbol{\mathcal{P}}_{v},\boldsymbol{\mathcal{M}})= \delta(\lim\limits_{k\rightarrow \infty}\boldsymbol{\mathcal{P}}_{v_{k}},\boldsymbol{\mathcal{P}}_{v_{2}})$,  $h_{i}=H(\mathcal{P}^{i}_{v_{1,i}},\mathcal{P}^{i}_{v_{2,i}})$  can be regarded as  the Hausdorff  distance between $ \mathcal{P}^{i}_{v_{i}}$ and $ \mathcal{M}^{i}$. Denote $ q_i=q_{2,i} $, then
$$
\delta(\boldsymbol{\mathcal{P}}_{v},\boldsymbol{\mathcal{M}})=\delta(\boldsymbol{H}(\boldsymbol{\mathcal{P}}_{v},\boldsymbol{\mathcal{M}}))
=r \sum_{i=1}^{N} \frac{q_{i}c_{i}}{\sqrt{\frac{2}{h_{i}\nu_{i}}-1}}.
$$
Finally, based on the definition of $ \epsilon$-GNE in Definition \ref{d3}, we analyze the difference  between $  J_{i}(\boldsymbol{x}^{*}(\boldsymbol{\mathcal{P}}_{v}))$ and $J_{i}(x_{i}^{\prime},\boldsymbol{x}^{*}_{-i}(\boldsymbol{\mathcal{P}}_{v})) $,  where the $ i $th player’s equilibrium strategy is $ x^{*}_{i}(\boldsymbol{\mathcal{P}}_{v}) $   with  respect to $ \boldsymbol{\mathcal{P}}_{v} $ and $ x^{\prime}_{i} $ is arbitrarily chosen from $ \mathcal{X}_{i} $. Meanwhile, other players’ strategies remain the same $ \boldsymbol{x}^{*}_{-i}(\boldsymbol{\mathcal{P}}_{v}) $.
$$
\small
\begin{aligned}
	&J_{i}\left(\boldsymbol{x}^{*}\left(\boldsymbol{\mathcal{P}}_{v}\right)\right)-J_{i}\left(x_{i}^{\prime}, \boldsymbol{x}_{-i}^{*}\left(\boldsymbol{\mathcal{P}}_{v}\right)\right) \\
	\leq &\left\|J_{i}\left(x_{i}^{\prime}, \boldsymbol{x}_{-i}^{*}(\boldsymbol{\mathcal{M}})\right)-J_{i}\left(x_{i}^{\prime}, \boldsymbol{x}_{-i}^{*}\left(\boldsymbol{\mathcal{P}}_{v}\right)\right)\right\|+\left\|J_{i}\left(\boldsymbol{x}^{*}\left(\boldsymbol{\mathcal{P}}_{v}\right)\right)-J_{i}\left(\boldsymbol{x}^{*}(\boldsymbol{\mathcal{M}})\right)\right\|\\
 &+J_{i}\left(\boldsymbol{x}^{*}(\boldsymbol{\mathcal{M}})\right)-J_{i}\left(x_{i}^{\prime}, \boldsymbol{x}_{-i}^{*}(\boldsymbol{\mathcal{M}})\right)\\
	\leq & \varsigma_{i}\left\|\boldsymbol{x}^{*}\left(\boldsymbol{\mathcal{P}}_{v}\right)-\boldsymbol{x}^{*}(\boldsymbol{\mathcal{M}})\right\|+\varsigma_{i}\left\|\boldsymbol{x}_{-i}^{*}(\boldsymbol{\mathcal{M}})-\boldsymbol{x}_{-i}^{*}\left(\boldsymbol{\mathcal{P}}_{v}\right)\right\|\\
	\leq& 2\varsigma_{i} \alpha_{1}^{-1}\left(\alpha_{2}\left(\alpha_{3}^{-1}\left(\frac{\delta\alpha_{4}(r)}{\mu}\right)\right)\right),
\end{aligned}
$$
where the third term in the first inequality is due to the definition of GNE. From this definition, the upper bound of the last term is zero.
This yields the conclusion.

\end{appendices}
\end{document}